\documentclass[12pt,a4paper,reqno]{article}    % Specifies the document style.
\usepackage{graphics}
\usepackage{cite}
\usepackage[demo]{graphicx}
\usepackage[fleqn]{amsmath}
\usepackage{booktabs}
\usepackage{multirow}
\usepackage{epsfig}
\usepackage{epstopdf}
\usepackage{booktabs}
\usepackage{placeins}
\usepackage{anysize}
\usepackage{enumerate}
\usepackage{latexsym,amssymb,amsmath,amscd,amsthm,amsxtra}
\usepackage{subfigure}
\usepackage{tabularx}
\usepackage[table]{xcolor}
\allowdisplaybreaks
\definecolor{lightgray}{gray}{0.9}
\newtheorem{theorem}{Theorem}

\newtheorem{exmp}{Example}

\numberwithin{equation}{section}

\def\({\left( }
\def\){\right )}
\begin{document}
\title {\textbf{An Efficient Computational Technique based on Cubic Trigonometric B-splines for  Time Fractional Burgers' Equation. }}
\author{Muhammad Yaseen, Muhammad Abbas\thanks{Corresponding authors: m.abbas@uos.edu.pk}\\
{Department of Mathematics, University of Sargodha, Pakistan.}\\
}
\date{}
\maketitle
\begin{abstract}
This paper presents a linear computational technique based on  cubic trigonometric cubic B-splines for time fractional burgers' equation. The nonlinear advection term is approximated by a new linearization technique which is very efficient and significantly reduces the computational cost. The usual finite difference formulation is used to approximate the Caputo time fractional derivative while the derivative in space is discretized using cubic trigonometric B-spline functions. The method is proved to be globally unconditionally stable. To measure the accuracy of the solution, a convergence analysis is also provided. A convergence analysis is Computational  experiments are performed  to further establish the accuracy and stability of the method. Numerical results are compared with those obtained by a scheme based on parametric spline functions. The comparison reveal that the proposed scheme is quite accurate and effective.  \\\\
\textbf{Key words:} Time fractional Burgers equation, Trigonometric basis functions, Cubic trigonometric B-splines method, Stability, Convergence.
\end{abstract}
\maketitle
%%%%%%%%%%%%%%%%%%%%%%%%%%%%%%%%%%%%%%%%%%%%%%%%%%%%%%%%%%%%%%%%%
% Type in your PAPER, starting below:
\section{Introduction}
We consider the following model of the nonlinnear time fractional Burgers' equation,
\begin{equation}\label{1.1}
    \frac{\partial^ \gamma u}{\partial t^\gamma}+ u \frac{\partial u}{\partial x}-\nu \frac{\partial^2 u}{\partial x^2}=0,\,\,\,\,\,\,a\leq x \leq b,\,\,\,\,\,0 < \gamma < 1,\,\,\,\,t>0.
\end{equation}
with initial condition
 \begin{equation}\label{1.2}
 u(x,0)=\phi(x)
\end{equation}
and the following boundary conditions
\begin{equation}\label{1.3}
\left\{
\begin{array}{cc}
 u(a,t)=\psi_1(t),\\
 u(b,t)=\psi_2(t)
\end{array}\qquad\qquad\qquad t>0,\right.
\end{equation}
where $a,b, \phi(x),\psi_1(t)$ and $\psi_2(t)$ are given, $\alpha>0$ is the reaction coefficient and  $\frac{\partial^\gamma}{\partial t^\gamma}u(x,t)$ represents the Caputo fractional derivative of order $\gamma$  given by \cite{pod,hil,kilbas}
\begin{equation}\label{1.4}
\frac{\partial^\gamma}{\partial t^\gamma}u(x,t)=
\begin{cases}
\frac{1}{\Gamma(1-\gamma)}\int\limits_0^t \frac{\partial u(x,s)}{\partial s}(t-s)^{-\gamma}ds,\qquad 0<\gamma <1\\
\frac{\partial u(x,s)}{\partial s},\qquad\qquad\qquad\qquad\qquad\qquad \gamma=1.
\end{cases}
\end{equation}
Burgers' equation ($\gamma=1$) is named for J. M. Burgers (1895-1981) and has numerous applications in various fields of science and engineering. As far as the fractional Burgers equation is concerned it describes the unidirectional propagation of weakly nonlinear acoustic waves through a pipe filled with gas. The fractional derivative appearing in the equation is responsible for the memory effects of the wall friction through the boundary layer. It is also used to model waves in bubbly liquids and shallow-water waves and consequently it is the main tool to model the shallow water problems. For more applications of the fractional Burgers' equation the reader is referred to \cite{garra,keller,djo,jero,sgimo}.

A number of numerical and analytical techniques are available in literature to solve fractional Burgers' equation. Adomian decomposition method has been utilized by Q. Wang \cite{wang}  to obtain approximate solutions of fractional Burgers' equation. Homotopy analysis method \cite{dhegan,song} and variational iteration method \cite{inc} have been applied on Burgers-type equations to find numerical solutions in the form of convergent series. A cubic parametric spline based numerical technique has been developed in \cite{denaf} for the solution of time fractional Burgers' equation. Cubic B-spline finite elements have been utilized in \cite{esen} to develop a numerical technique for the solutions of time fractional Burgers' equation. Monami \cite{momani} obtained non-perturbative analytic solutions of the space- and time-fractional Burgers' equation. Hereditary effects on nonlinear acoustic waves have been discussed by Sugimoto \cite{sg} using different models of Burgers' equation with a fractional derivative. A generalized differential transform method has been used in \cite{liu} to obtain numerical solutions of the space- and time-fractional coupled Burgers' equations. Xhu and Zhang \cite{liu1} have presented a new finite difference scheme for generalized Rosenau-Burgers equation. Scaling limit solution of a fractional Burgers' equation has been obtained in \cite{liu2}.

Finite difference schemes are most widely used for solving a variety of fractional differential equations because of simplicity and perception \cite{zhu1,zhu2,zhu3,zhu4,zhu5,zhu6}. Stability and convergence analysis of various finite difference schemes have been discussed in \cite{zhu7,zhu8}. It has been pointed out there that in explicit difference schemes small step sizes should be chosen for keeping the method stable. In case of implicit schemes, the stability may be unconditional but one obtains a new system of nonlinear algebraic equations  at each time level. The main aim of this paper is to develop a numerical scheme which is computationally efficient and  surmount the above mentioned difficulties for time fractional Burgers' equation. The presented scheme is shown to be unconditionally globally stable. A convergence analysis of the scheme is also presented to ensure the accuracy of the scheme. Compared to the existing numerical schemes, computational cost is greatly reduced by discretizing the nonlinear advection term by a new technique. Numerical experiments are carried out to further establish the accuracy and validity of the scheme.

The remainder of the paper is organized as follows. In section 2, the numerical scheme based on cubic trigonometric B-splines is derived in detail. Section 3 discusses the stability and convergence of the scheme is discussed in section 4. Section 5 shows comparison of our numerical results with those of \cite{denaf}. Section 5 summarizes the conclusions of this study.

\section{The Derivation of the Scheme}
For given positive integers $M$ and $N$, let $\tau=\frac{T}{N}$ be the temporal and $h=\frac{L}{M}$ spatial step sizes respectively. Following the usual notations, set $t_n=n\tau~(0\leq n \leq N)$, $x_j=jh,~(0 \leq j \leq M)$, $\Omega_{\tau}=\{t_n|0 \leq n \leq N\}$ and $\Omega_h=\{x_j|0 \leq n \leq N\}$. Let $u_j^n$ be approximation to exact solution at the point $(x_j,t_n)$ and $\mathcal{A}=\{u_j^n|0\leq j\leq M, 0\leq n\leq N\}$ be grid function space defined on $\Omega_h \times \Omega_h$.  The solution domain $0\leq x\leq L$ is uniformly partitioned  by knots $x_{i}$ into $N$ subintervals $[x_{i}, x_{i+1}]$ of equal length $h$, $i=0,1,2,...,N-1$, where $0=x_{0}<x_{1}<...<x_{n-1}<x_{N}=L$. Our scheme for solving (\ref{1.1}) requires approximate solution $U(x,t)$ to the exact solution $u(x,t)$  in the following form  \cite{36,37}
\begin{equation}\label{3.1}
    U(x,t)=\sum\limits_{i=-1}^{N-1} c_i(t) TB_i^4 (x),
\end{equation}
where $c_i(t)$ are unknowns to be determined and $TB^4_i (x)$  \cite{abbas} are twice differentiable cubic Trigonometric basis functions given by
\begin{equation}\label{3.2}
    TB^4_i (x)=\frac{1}{w}
    \begin{cases}
    p^3(x_i) & x\in [x_i,x_{i+1}]\\
    p(x_i)(p(x_i)q(x_{i+2})+q(x_{i+3})p(x_{i+1}))+q(x_{i+4})p^2(x_{i+1}),  & x \in [x_{i+1},x_{i+2}]\\
    q(x_{i+4})(p(x_{i+1})q(x_{i+3})+q(x_{i+4})p(x_{i+2}))+p(x_i)q^2(x_{i+3}), & x \in [x_{i+2},x_{i+3}]\\
    q^3(x_{i+4}), & x \in [x_{i+3},x_{i+4}]
    \end{cases}
\end{equation}
where\\
$\displaystyle{p\(x_i\)=\sin\(\frac{x-x_i}{2}\), q\(x_i\)=\sin\(\frac{x_i-x}{2}\), w=\sin\(\frac{h}{2}\) \sin\(h\) \sin\(\frac{3 h}{2}\).}$\\
 Due to local support property of the cubic trigonometric B-splines only $TB^4_{j-1}(x), TB^4_{j}(x)$ and $TB^4_{j+1}(x)$ are survived so that the approximation $u_j^n$ at the grid point $(x_j,t_n)$ at $n$ th time level is given as:
\begin{equation}\label{3.3}
    u(x_{j},t_{n})=u_j^n=\sum\limits_{j=i-1}^{i+1} c_j^n(t) TB^4_j (x).
\end{equation}
The  time dependent unknowns $c_j^n(t)$ are to be determined by making use of the initial and boundary conditions, and the collocation conditions on $TB^4_i (x)$. As a result the approximations $u_j^n$ and its necessary derivatives as given below:
\begin{equation}\label{1.8}
    \begin{cases}
   \displaystyle{ u_j^n=a_1 c_{j-1}^n+a_2 c_{j}^n+ a_1 c_{j+1}^n}\\
    \displaystyle{ (u_j^n)_x=-a_3 c_{j-1}^n+a_3 c_{j+1}^n}\\
    \displaystyle{ (u_j^n)_{xx}=a_4 c_{j-1}^n+ a_5 c_{j}^n +a_4 c_{j+1}^n,}
    \end{cases}
\end{equation}
where\\
$\displaystyle{a_1=\csc\(h\) \csc\(\frac{3h}{2}\)\sin^2\(\frac{h}{2}\),}$\\
$\displaystyle{a_2=\frac{2}{1+2 \cos\(h\)},}$\\
$\displaystyle{a_3=\frac{3}{4} \csc \(\frac{3h}{2}\),}$\\
$\displaystyle{a_4=\frac{3+9\cos\(h\)}{4 \cos\(\frac{h}{2}\)-4 \cos\(\frac{5 h}{2}\)},}$\\
$\displaystyle{a_5=-\frac{3 \cot^2\(\frac{h}{2}\)}{2+4 \cos\(h\)}.}$\\\\
Note that the nonlinear advection term $u\frac{\partial u}{\partial x}$ is usually discretized as $u_j^n (u_x)_j^{n+1}+u_j^{n+1}(u_x)_j^n-u_j^n(u_x)_j^n$. Instead, we present a new technique to linearize the advection term as follows:
\begin{equation}\label{1.9}
    u\frac{\partial u}{\partial x}|_{(x_j,t_n)}=\frac{1}{3}\left[u_j^n(u_x)_j^{n+1}+(u_j^nu_j^{n+1})_x\right]
\end{equation}
which is simplified to
\begin{equation}\label{1.10}
     u\frac{\partial u}{\partial x}|_{(x_j,t_n)}=\frac{1}{3}u_j^{n+1}(u_x)_j^n+\frac{2}{3}u_j^n(u_x)_j^{n+1}
\end{equation}
Following \cite{sidd}, the fractional derivative  $\frac{\partial^\gamma}{\partial t^\gamma}u(x,t)$ is discrtetized by
\begin{equation}\label{1.11}
    \frac{\partial^\gamma}{\partial t^\gamma}u(x,t)|_{(x_j,t_{n+1})}=\frac{\partial^\gamma}{\partial t^\gamma}u_j^{n+1}=\frac{1}{\Delta t^\gamma \Gamma[2-\gamma]}\sum\limits_{l=0}^{n}w_l(u_j^{n-l+1}-u_j^{n-l}),
\end{equation}
where $w_l=(l+1)^{1-\gamma}-l^{1-\gamma}$. It is straight forward to confirm that
\begin{itemize}
  \item $w_l>0,~~l=0,1,\cdots ,n$,
  \item $1=w_0>w_1>w_2>\cdots,w_n~\text{and}~w_n\rightarrow 0~\text{as}~n\rightarrow \infty$.
  \item $\sum\limits_{l=0}^n (w_l-w_{l+1})=1$.
\end{itemize}
After some simplifications, equation (\ref{1.11}) can be re written as
\begin{equation}\label{1.12}
    \frac{\partial^\gamma}{\partial t^\gamma}u(x_j,t_{n+1})=\frac{1}{\alpha_0}\left[u_j^{n+1}+\sum\limits_{l=0}^{n} (w_{n-l+1}-w_{n-l})u_j^l - w_n u_j^0\right],
\end{equation}
where $\alpha_0=\Delta t^\gamma \Gamma[2-\gamma]$.
Using (\ref{1.12}), the time fractional burgers' equation (\ref{1.1}) can be written in discretized form as:
\begin{equation}\label{1.13}
    u_j^{n+1}=\sum\limits_{l=0}^{n} (w_{n-l}-w_{n-l+1})u_j^l + w_n u_j^0+\nu\alpha_0 (u_j^{n+1})_{xx}-\frac{\alpha_0}{3}\left( u_j^n(u_j^{n+1})_x+(u_j^nu_j^{n+1})_x\right).
\end{equation}
Now we introduce some important notations necessary to understand the subsequent material. For any grid function $u \in  \mathcal{A}$, define\\
$(u^n,v^n)=h\sum\limits_{j=1}^{M-1} u_j^nv_j^n,~~~\|u^n\|^2=(u^n,u^n),~~~\|u^n\|_\infty=\max\limits_{0\leq j\leq M-1} |u_j^n|,$\\$((u^n)_{xx},u^n)=-((u^n)_x,(u^n)_x)=-\|(u^n)_x\|^2.$
\section{Stability analysis}
In this section, the scheme (\ref{1.13}) is proved to be unconditionally globally stable.
\begin{theorem}\label{th1}
The scheme (\ref{1.13}) is bounded that is there exists a constant $C$ such that $$\|u_j^{n+1}\|^2\leq C,~~~n=0,1,\cdots,N.$$
\end{theorem}
\begin{theorem}
The scheme (\ref{1.13}) is globally unconditionally stable.
\end{theorem}
%\noindent
%\textbf{Proof:} The conclusion follows directly from the boundedness of the numerical scheme (\ref{1.13}).
\section{Convergence Analysis}
Assuming $H>0$  which assumes different values at various locations and is independent of $j, n, h,~\text{and}~\tau$ with the condition that
 \begin{equation}\label{4.1}
    |u_{tt}|\leq H,~|u_{xxxx}|\leq H ~\text{for all}~ (x,t)\in \Omega_h\times \Omega_\tau.
\end{equation}
 Then for the scheme (\ref{1.13}), we have
\begin{equation}\label{4.2}
\begin{split}
    u(x_j,t^{n+1})=\sum\limits_{l=0}^{n} (w_{n-l}-w_{n-l+1})u(x_j,t^{l}) + w_n u(x_j,t^{0})+\nu\alpha_0 (u(x_j,t^{n+1}))_{xx}-\\{}\frac{\alpha_0}{3}\left( u(x_j,t^{n})(u(x_j,t^{n+1}))_x+(u(x_j,t^{n})u(x_j,t^{n+1}))_x\right),
    \end{split}
\end{equation}
where $(x_j,t^{n})$ is the exact solution at point $(x_j,t^{n})$.
\begin{theorem}
If $u(x,t)$ satisfies smoothness condition (\ref{4.1}), then for sufficiently small $h$ and $\tau$,
\begin{equation}\label{4.3}
    \|e^{n+1}\|\leq D+O(\tau^2+\tau h^2),
\end{equation}
where $e_j^{n+1}=u(x_j,t^{n+1})-u_j^{n+1}~~\text{and}~~D=\max\limits_{0\leq l \leq n} \|e^l\|$
\end{theorem}
\section{Numerical experiments and discussion}
For the purpose of numerical computations, we plug in values from equation (\ref{1.8}) into equation (\ref{1.13}). As a result the scheme (\ref{1.13}) is written equivalently as:
\begin{equation}\label{4.1}
         \begin{split}
       L_j^n c_{j-1}^{n+1}+M_j^n c_{j}^{n+1}+N_j^n c_{j+1}^{n+1}=\sum\limits_{l=1}^{n} (b_{n-l}-b_{n+1-l})(a_1c_{j-1}^{l}+a_2c_{j}^{l}+a_1c_{j+1}^{l})\\{}+b_n(a_1c_{j-1}^{0}+a_2c_{j}^{0}+a_1c_{j+1}^{0}),
      \end{split}
     \end{equation}
     where,
     $$L_j^n=a_1+\frac{a_1\alpha_0}{3}(-b_1c_{j-1}^{n}+b_2c_{j}^{n}+b_1c_{j+1}^{n})-\frac{2b_1\alpha_0}{3}(a_1c_{j-1}^{n}+a_2c_{j}^{n}+a_1c_{j+1}^{n})-\nu\alpha_0c_1,$$
     $$M_j^n=a_2+\frac{a_2\alpha_0}{3}(-b_1c_{j-1}^{n}+b_2c_{j}^{n}+b_1c_{j+1}^{n})+\frac{2b_2\alpha_0}{3}(a_1c_{j-1}^{n}+a_2c_{j}^{n}+a_1c_{j+1}^{n})-\nu\alpha_0c_2,$$
     and
     $$N_j^n=a_1+\frac{a_1\alpha_0}{3}(-b_1c_{j-1}^{n}+b_2c_{j}^{n}+b_1c_{j+1}^{n})+\frac{2b_1\alpha_0}{3}(a_1c_{j-1}^{n}+a_2c_{j}^{n}+a_1c_{j+1}^{n})-\nu\alpha_0c_1.$$
The linear system (\ref{3.1}) consists of  contains $(N-1)$ equations in $(N+1)$ unknowns. To obtain a unique solution to the system, we need two additional equations which can be obtained by utilizing the given  boundary conditions (\ref{1.3}). As a result,  a  matrix system of dimension $(N + 1)\times(N + 1)$  is obtained which is a tri-diagonal  system which can be uniquely solved by Thomas Algorithm \cite{38}.
\begin{exmp}
Consider the Burgers' equation \cite{denaf}
\begin{equation}\label{3.2}
     \frac{\partial^ \gamma u}{\partial t^\gamma}+ u \frac{\partial u}{\partial x}-\nu \frac{\partial^2 u}{\partial x^2}=0,\,\,\,\,\,\,-3\leq x \leq 3,\,\,\,\,\,0 < \gamma < 1,\,\,\,\,t>0.
\end{equation}
subject to the boundary condition
\begin{equation}\label{3.3}
    u(a,t)=\theta_1(t),~~u(b,t)=\theta_2(t)
\end{equation}
and initial condition
\begin{equation}\label{3.4}
    u(x,0)=\frac{\mu+\sigma+(\sigma-\mu)\exp[\frac{\mu}{\nu}(x-\lambda)]}{1+\exp[\frac{\mu}{\nu}(x-\lambda)]}
\end{equation}
\end{exmp}
The exact solution of this problem is
\begin{equation}\label{3.5}
    u(x,t)=\frac{\mu+\sigma+(\sigma-\mu)\exp[\frac{\mu}{\nu}(x-\sigma t-\lambda)]}{1+\exp[\frac{\mu}{\nu}(x-\sigma t-\lambda)]}
\end{equation}
The approximate solution is obtained by applying the presented scheme (\ref{4.1}) and the results are presented in Tables 1-3 with $\mu=0.3, ~\sigma=0.4,~ \nu=0.1,~ \lambda=0.8,~h=0.01=\tau$. Figure 1 plots the behavior of approximate solutions at different time levels. The graphs are in excellent agreement with those of \cite{denaf}. In figure 2, a 3-D comparison between the exact and approximate solutions shows the tremendous accuracy of the scheme.  Figure 3 depicts the absolute error profile at $t=1$.  In Tables 1-3, a comparative study of the obtained results with those of \cite{denaf} is presented. The comparison reveal that the presented scheme provides better accuracy and efficiency.

\begin{exmp}
Consider the fractional Burgers' equation \cite{denaf}
\begin{equation}\label{4.6}
\frac{\partial^ \gamma u}{\partial t^\gamma}+ u \frac{\partial u}{\partial x}-\nu \frac{\partial^2 u}{\partial x^2}=0,\,\,\,\,\,\,a\leq b \leq 3,\,\,\,\,\,0 < \gamma < 1,\,\,\,\,t>0.
\end{equation}
with initial condition
\begin{equation}\label{4.7}
     u(x,0)=\frac{\mu+\sigma+(\sigma-\mu)\exp[\frac{\mu}{\nu}(x-\lambda)]}{1+\exp[\frac{\mu}{\nu}(x-\lambda)]}
\end{equation}
and the boundary conditions
\begin{equation}\label{4.8}
    u(-3,t)\approx 0.699993+(1.07 \times 10^{-5})\frac{t^{\gamma}}{\Gamma(1+\gamma)}-(9.67 \times 10^{-6})\frac{t^{2\gamma}}{\Gamma(1+2\gamma)}+(1.16 \times 10^{-5})\frac{t^{3\gamma}}{\Gamma(1+3\gamma)}
\end{equation}
\begin{equation}\label{4.9}
     u(3,t)\approx 0.100815+(1.3 \times 10^{-3})\frac{t^{\gamma}}{\Gamma(1+\gamma)}+(1.17 \times 10^{-3})\frac{t^{2\gamma}}{\Gamma(1+2\gamma)}-(5.72 \times 10^{-6})\frac{t^{3\gamma}}{\Gamma(1+3\gamma)}
\end{equation}
\end{exmp}
The presented scheme (\ref{4.1}) is employed to approximate the solution when $\mu=0.3,~ \sigma=0.4,~\lambda=0.8,~\nu=0.1$ and the obtained results are presented in Table 4 and Figure 4. In Table 4 approximate solutions are compared with those obtained in \cite{denaf} at time $\gamma=0.2$ and $\gamma=0.8$. The results are in close agreement. Figure 4 depicts the behavior of numerical solution when $\gamma=0.8,~\gamma=0.5$ and $\gamma=0.2$. The graphs are in tremendous agreement with those of \cite{denaf}.
\section{Concluding Remarks}\label{s6}
This study presents a linear numerical technique based on cubic Trigonometric B-splines for the time fractional Burgers' equation. The algorithm utilizes usual finite difference schemes to approximate the Caputo time fractional derivative and the derivative in space are approximated using the cubic trigonometric B-spline basis functions. A new approximation techniques is presented for nonlinear advection term which is highly efficient and tremendously reduces amount of numerical calculations. A special attention has been given to study the stability and convergence of the scheme. The scheme is shown to be globally unconditionally stable. The obtained results are compared with those obtained by a numerical techniques based on parametric spline functions. Numerical and graphical comparisons reveal that the presented scheme is better and computationally very efficient. Consequently, the scheme can be applied to similar type of nonlinear fractional partial differential equations.

%%%%%%%%%%%%%%%
\end{document}